\documentclass[review,english]{elsarticle_mod}

\usepackage{lineno,hyperref}
\modulolinenumbers[5]
\PassOptionsToPackage{normalem}{ulem}
\usepackage{ulem}
\usepackage{babel}
\usepackage{color}
\usepackage{amsmath}
\usepackage{amssymb}
\usepackage{amsthm}
\usepackage{geometry}
\makeatletter
\numberwithin{equation}{section}
\numberwithin{figure}{section}
  \theoremstyle{remark}
  \newtheorem*{notation*}{\protect\notationname}
\theoremstyle{plain}
\newtheorem{thm}{\protect\theoremname}
  \theoremstyle{definition}
  
  \theoremstyle{plain}
  
  \theoremstyle{remark}
  \newtheorem{rem}[thm]{\protect\remarkname}

\makeatother

  \providecommand{\definitionname}{Definition}
  \providecommand{\lemmaname}{Lemma}
  \providecommand{\notationname}{Notation}
  \providecommand{\remarkname}{Remark}
\providecommand{\theoremname}{Theorem}

\begin{document}

\begin{frontmatter}

\title{A Note on Iterations-based Derivations of High-order Homogenization Correctors for Multiscale Semi-linear Elliptic Equations}


\author[aff1]{Vo Anh Khoa\corref{mycorrespondingauthor}}
\cortext[mycorrespondingauthor]{Corresponding author}
\ead{khoa.vo@gssi.infn.it, vakhoa.hcmus@gmail.com}

\author[aff2]{Adrian Muntean}
\ead{adrian.muntean@kau.se}

\address[aff1]{Mathematics and Computer Science Division,
Gran Sasso Science Institute, L'Aquila, Italy}
\address[aff2]{Department of Mathematics and Computer Science, Karlstad University,
Sweden}

\begin{abstract}
This Note aims at presenting a simple and efficient procedure to derive the structure of high-order corrector estimates for the homogenization limit  applied to a semi-linear elliptic equation posed in perforated domains. Our working technique relies on monotone iterations combined with formal two-scale homogenization asymptotics. It can be adapted to handle more complex scenarios  including for instance nonlinearities posed  at the boundary of perforations and  the vectorial case, when the model equations are coupled only through the nonlinear production terms.
\end{abstract}


\end{frontmatter}

\linenumbers

\section{Background}

Modern approaches to modeling focus on multiple scales. Given a multiscale physical  problem, one of the leading  questions is to derive upscaled model equations and the corresponding structure of effective model coefficients (e.g. \cite{Wei11,SPK14}).   This Note aims at exploring the quality of the upscaling/homogenization procedure by deriving whenever possible corrector (error) estimates for the involved unknown functions, fields, etc.  and their gradients (i.e. of the transport fluxes).  Ultimately, these estimates contribute essentially to the control of the approximation error of numerical methods to  multiscale PDE problems. 

Our starting point is a microscopic PDE model describing the motion of populations of colloidal particles in soils and  porous tissues with direct applications in drug-delivery design and control of the spread of radioactive pollutants (\cite{ray2013thesis,AWR,RMK12}).   We have previously analyzed a reduced version of this system in \cite{KM15}. Here,  we point out a short  alternative proof based on monotone iterations (\cite{Pao93}) of the corrector estimates derived in \cite{KM15} and extend them to higher asymptotic orders.

\section{Problem setting}

We are concerned with the study of the semi-linear elliptic boundary
value problem of the form
\begin{equation}
\begin{cases}
\mathcal{A}^{\varepsilon}u^{\varepsilon}=R\left(u^{\varepsilon}\right), & x\in\Omega^{\varepsilon},\\
u^{\varepsilon}=0, & x\in\Gamma^{ext},\\
\nabla u^{\varepsilon}\cdot\mbox{n}=0, & x\in\Gamma^{\varepsilon},
\end{cases}\label{eq:1.1}
\end{equation}
where the operator $\mathcal{A}^{\varepsilon}u^{\varepsilon}:=\nabla\cdot\left(-d^{\varepsilon}\nabla u^{\varepsilon}\right)$
involves $d^{\varepsilon}$ termed as the molecular diffusion while
$R$ represents the volume reaction rate. We take into account 
the following assumptions:

$\left(\mbox{A}_{1}\right)$ the diffusion coefficient $d^{\varepsilon}\in L^{\infty}\left(\mathbb{R}^{d}\right)$
for $d\in\mathbb{N}$ is $Y$-periodic and symmetric, and it guarantees
the ellipticity of $\mathcal{A}^{\varepsilon}$ as follows:
\[
d^{\varepsilon}\xi_{i}\xi_{j}\ge\alpha\left|\xi\right|^{2}\quad\mbox{for any}\;\xi\in\mathbb{R}^{d};
\]

$\left(\mbox{A}_{2}\right)$ the reaction coefficient $R\in L^{\infty}\left(\Omega^{\varepsilon}\times\mathbb{R}\right)$
is globally $L-$Lipschitzian, i.e. there exists $L>0$ independent
of $\varepsilon$ such that
\[
\left\Vert R\left(u\right)-R\left(v\right)\right\Vert \le L\left\Vert u-v\right\Vert \quad\mbox{for}\;u,v\in\mathbb{R}.
\]

\begin{figure}[!h]
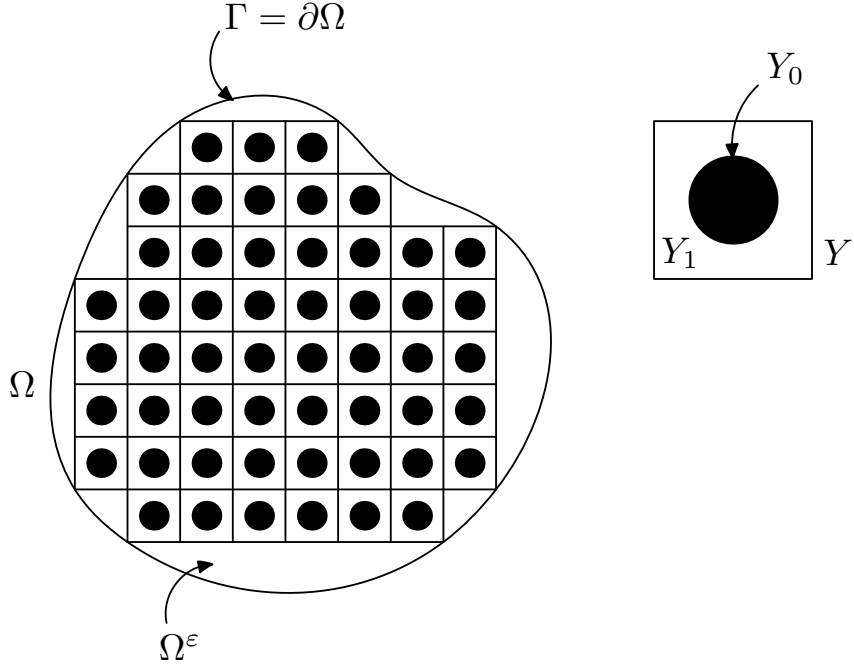
 	
\centering 	
	\parbox{10cm}{\convertMPtoPDF{fig.man}{1.4}{1.4}}	
	\caption{An admissible perforated domain (left) and basic geometry of the microstructure (right).}
	\label{fig:1}
\end{figure}

It is worth noting that the domain $\Omega^{\varepsilon}\in\mathbb{R}^{d}$
considered here approximates a porous medium.  The precise description  of $\Omega^{\varepsilon}$ is
showed in \cite[Section 2]{KM15} and \cite{HJ91}. In
Figure \ref{fig:1} (left), we sketch an admissible geometry of our
medium, pointing out the sample microstructure in Figure \ref{fig:1} (right). We follow the notation from \cite{KM15}.

\begin{rem}
In this paper, we denote the space $V^{\varepsilon}$ by
\[
V^{\varepsilon}:=\left\{v\in H^{1}(\Omega^{\varepsilon})|v=0\;\mbox{on}\;\Gamma^{ext}\right\}
\]
endowed with the norm
\[
\left\Vert v\right\Vert _{V^{\varepsilon}}=\left(\int_{\Omega^{\varepsilon}}\left|\nabla v\right|^{2}dx\right)^{1/2}.
\]
This norm is equivalent (uniformly in the homogenization parameter $\varepsilon$) to the usual $H^1$ norm by the Poincar\'e inequality.
\end{rem}

\section{Derivation of corrector estimates}

We introduce the following $M$th-order expansion $\left(M\ge2\right)$:
\begin{equation}
u^{\varepsilon}\left(x\right)=\sum_{m=0}^{M}\varepsilon^{m}u_{m}\left(x,y\right)+\mathcal{O}\left(\varepsilon^{M+1}\right),\quad x\in\Omega^{\varepsilon},\label{eq:expansion}
\end{equation}
where $u_{m}\left(x,\cdot\right)$ is $Y$-periodic for $0\le m\le M$.

Following standard homogenization procedures, we  deduce the so called auxiliary problems (see e.g. \cite{BLP11}). 
To do so, we consider the functional $\Phi\left(x,y\right)$ depending
on two variables: the macroscopic $x$ and $y=x/\varepsilon$ the
microscopic presentation, and denote by $\Phi^{\varepsilon}\left(x\right)=\Phi\left(x,y\right)$.
The simple chain rule allows us to derive the fact that
\begin{equation}
\nabla\Phi^{\varepsilon}\left(x\right)=\nabla_{x}\Phi\left(x,\frac{x}{\varepsilon}\right)+\varepsilon^{-1}\nabla_{y}\Phi\left(x,\frac{x}{\varepsilon}\right).\label{eq:chain}
\end{equation} The quantities $\nabla u^{\varepsilon}$ and $\mathcal{A}^{\varepsilon}u^{\varepsilon}$
must be expended correspondingly. In fact, it follows from \eqref{eq:chain} and \eqref{eq:expansion}
that
\begin{eqnarray}
\nabla u^{\varepsilon} & = & \left(\nabla_{x}+\varepsilon^{-1}\nabla_{y}\right)\left(\sum_{m=0}^{M}\varepsilon^{m}u_{m}+\mathcal{O}\left(\varepsilon^{M+1}\right)\right)\nonumber \\
 & = & \varepsilon^{-1}\nabla_{y}u_{0}+\sum_{m=0}^{M-1}\varepsilon^{m}\left(\nabla_{x}u_{m}+\nabla_{y}u_{m+1}\right)+\mathcal{O}\left(\varepsilon^{M}\right).\label{eq:laplaceu}
\end{eqnarray}

Using the structure of the operator $\mathcal{A}^{\varepsilon}$, we
obtain the following:
\begin{eqnarray}
\mathcal{A}^{\varepsilon}u^{\varepsilon} & = & \varepsilon^{-2}\nabla_{y}\cdot\left(-d\left(y\right)\nabla_{y}u_{0}\right)\nonumber \\
 &  & +\varepsilon^{-1}\left[\nabla_{x}\cdot\left(-d\left(y\right)\nabla_{y}u_{0}\right)+\nabla_{y}\cdot\left(-d\left(y\right)\left(\nabla_{x}u_{0}+\nabla_{y}u_{1}\right)\right)\right]\nonumber \\
 &  & +\sum_{m=0}^{M-2}\varepsilon^{m}\left[\nabla_{x}\cdot\left(-d\left(y\right)\left(\nabla_{x}u_{m}+\nabla_{y}u_{m+1}\right)\right)\right.\nonumber \\
 &  & \left.+\nabla_{y}\cdot\left(-d\left(y\right)\left(\nabla_{x}u_{m+1}+\nabla_{y}u_{m+2}\right)\right)\right]+\mathcal{O}\left(\varepsilon^{M-1}\right).\label{eq:newoperator}
\end{eqnarray}

Concerning the boundary condition at $\Gamma^{\varepsilon}$, we note:
\begin{equation}
d^{\varepsilon}\nabla u^{\varepsilon}\cdot\mbox{n}=d_{i}\left(y\right)\left(\varepsilon^{-1}\nabla_{y}u_{0}+\sum_{m=0}^{M-1}\varepsilon^{m}\left(\nabla_{x}u_{m}+\nabla_{y}u_{m+1}\right)\right)\cdot\mbox{n}.\label{eq:boundary}
\end{equation}

To investigate the convergence analysis, we consider the following structural property:
\begin{equation}
R\left(\sum_{m=0}^{M}\varepsilon^{m}u_{m}\right)=\sum_{m=0}^{M}\varepsilon^{m-r}R\left(u_{m}\right)+\mathcal{O}\left(\varepsilon^{M-r+1}\right)\quad\mbox{for}\;r\in\mathbb{Z},r\le2.\label{eq:newR}
\end{equation}

At this point we see, if $r\in\left\{ 1,2\right\} $ it then generate
spontaneously nonlinear auxiliary problems. To see the impediment,
let us focus on $r=2$. By collecting the coefficients of the same
powers of $\varepsilon$ in \eqref{eq:newoperator} and \eqref{eq:boundary},
we are led to the following systems, which we refer to the auxiliary
problems:
\begin{equation}
\begin{cases}
\mathcal{A}_{0}u_{0}=R\left(u_{0}\right), & \mbox{in}\;Y_{1},\\
-d\left(y\right)\nabla_{y}u_{0}\cdot\mbox{n}=0, & \mbox{on}\;\partial Y_{0},\\
u_{0}\;\mbox{is}\;Y-\mbox{periodic in}\;y,
\end{cases}\label{eq:prob1}
\end{equation}

\begin{equation}
\begin{cases}
\mathcal{A}_{0}u_{1}=R\left(u_{1}\right)-\mathcal{A}_{1}u_{0}, & \mbox{in}\;Y_{1},\\
-d\left(y\right)\left(\nabla_{x}u_{0}+\nabla_{y}u_{1}\right)\cdot\mbox{n}=0, & \mbox{on}\;\partial Y_{0},\\
u_{1}\;\mbox{is}\;Y-\mbox{periodic in}\;y,
\end{cases}\label{eq:prob2}
\end{equation}

\begin{equation}
\begin{cases}
\mathcal{A}_{0}u_{m+2}=R\left(u_{m+2}\right)-\mathcal{A}_{1}u_{m+1}-\mathcal{A}_{2}u_{m}, & \mbox{in}\;Y_{1},\\
-d\left(y\right)\left(\nabla_{x}u_{m+1}+\nabla_{y}u_{m+2}\right)\cdot\mbox{n}=0, & \mbox{on}\;\partial Y_{0},\\
u_{m+2}\;\mbox{is}\;Y-\mbox{periodic in}\;y,
\end{cases}\label{eq:probm}
\end{equation}
for $0\le m\le M-2$.

Here, we have denoted by
\begin{eqnarray*}
\mathcal{A}_{0} & := & \nabla_{y}\cdot\left(-d\left(y\right)\nabla_{y}\right),\\
\mathcal{A}_{1} & := & \nabla_{x}\cdot\left(-d\left(y\right)\nabla_{y}\right)+\nabla_{y}\cdot\left(-d\left(y\right)\nabla_{x}\right),\\
\mathcal{A}_{2} & := & \nabla_{x}\cdot\left(-d\left(y\right)\nabla_{x}\right).
\end{eqnarray*}

\begin{rem}
In the case $r\le0$, it is trivial to not only prove the well-posedness
of these auxiliary problems \eqref{eq:prob1}-\eqref{eq:probm}, but
also to compute the solutions by many approaches due to its linearity.
For details, the reader is referred here  to \cite{CP99}.
\end{rem}
The idea is now to linearize the auxiliary problems.  Inspired by the fact that a fixed-point homogenization argument seems to be applicable in this framework, and also by the way a priori  error estimates are proven for difference schemes, we suggest an iteration technique to "linearize" the involved PDE systems.  We start the procedure by choosing  the initial
point $u_{m}^{\left(0\right)}=0$ for $m\in\left\{ 0,...,M\right\}$.
As next step, we consider the following systems corresponding to the nonlinear
auxiliary problems:
\begin{equation}
\begin{cases}
\mathcal{A}_{0}u_{0}^{\left(n_{0}\right)}=R\left(u_{0}^{\left(n_{0}-1\right)}\right), & \mbox{in}\;Y_{1},\\
-d\left(y\right)\nabla_{y}u_{0}^{\left(n_{0}\right)}\cdot\mbox{n}=0, & \mbox{on}\;\partial Y_{0},\\
u_{0}^{\left(n_{0}\right)}\;\mbox{is}\;Y-\mbox{periodic in}\;y,
\end{cases}\label{eq:prob1-1}
\end{equation}

\begin{equation}
\begin{cases}
\mathcal{A}_{0}u_{1}^{\left(n_{1}\right)}=R\left(u_{1}^{\left(n_{1}-1\right)}\right)-\mathcal{A}_{1}u_{0}^{\left(n_{0}\right)}, & \mbox{in}\;Y_{1},\\
-d\left(y\right)\left(\nabla_{x}u_{0}^{\left(n_{0}\right)}+\nabla_{y}u_{1}^{\left(n_{1}\right)}\right)\cdot\mbox{n}=0, & \mbox{on}\;\partial Y_{0},\\
u_{1}^{\left(n_{1}\right)}\;\mbox{is}\;Y-\mbox{periodic in}\;y,
\end{cases}\label{eq:prob2-1}
\end{equation}

\begin{equation}
\begin{cases}
\mathcal{A}_{0}u_{m+2}^{\left(n_{m+2}\right)}=R\left(u_{m+2}^{\left(n_{m+2}-1\right)}\right)-\mathcal{A}_{1}u_{m+1}^{\left(n_{m+1}\right)}-\mathcal{A}_{2}u_{m}^{\left(n_{m}\right)}, & \mbox{in}\;Y_{1},\\
-d\left(y\right)\left(\nabla_{x}u_{m+1}^{\left(n_{m+1}\right)}+\nabla_{y}u_{m+2}^{\left(n_{m+2}\right)}\right)\cdot\mbox{n}=0, & \mbox{on}\;\partial Y_{0},\\
u_{m+2}^{\left(n_{m+2}\right)}\;\mbox{is}\;Y-\mbox{periodic in}\;y,
\end{cases}\label{eq:probm-1}
\end{equation}
for $0\le m\le M-2$. Note that the quantity $n_{m}$ is independent
of $\varepsilon$.

Since the approximate auxiliary problems became linear, standard procedures are able to find the solutions $u_{m}^{\left(n_{m}\right)}$
for $0\le m\le M$ . Note that these
problems admit a unique solution (see, e.g. \cite[Lemma 2.2]{CP99})
on $V$, i.e. the quotient space of $V_{Y_{1}}$ defined by
\[
V_{Y_{1}}:=\left\{ \varphi|\varphi\in H^{1}\left(Y_{1}\right),\varphi\;\mbox{is}\;Y-\mbox{periodic}\right\} .
\]

On the other side, if $\kappa_{p}:=C_{p}L\alpha^{-1}<1$ holds (here  $C_{p}$ is the Poincar\'e constant depending only on the dimension
of $Y_{1}$), we easily obtain that for every $m$, $\left\{ u_{m}^{\left(n_{m}\right)}\right\} $
is a Cauchy sequence in $H^{1}\left(Y_{1}\right)$. Hereby, it naturally
claims the existence and uniqueness of the nonlinear auxiliary problems
\eqref{eq:prob1}-\eqref{eq:probm}. Moreover, the convergence rate
of the iteration procedure is given by
\[
\left\Vert u_{m}^{\left(n_{m}\right)}-u_{m}\right\Vert _{H^{1}\left(Y_{1}\right)}\le\frac{\kappa_{p}^{n_{m}}}{1-\kappa_{p}^{n_{m}}}\left\Vert u_{m}^{\left(1\right)}\right\Vert _{H^{1}\left(Y_{1}\right)}.
\]

\begin{rem}
For more details in this sense, see \cite[Theorem 9]{KM15}
and \cite[Theorem 2.2]{KTLN15}.
\end{rem}
To prove the corrector estimate, we suppose that the solutions
of the auxiliary problems \eqref{eq:prob1}-\eqref{eq:probm} belong
to the space $L^{\infty}\left(\Omega^{\varepsilon};V\right)$. Let
us introduce the following function:
\[
\varphi^{\varepsilon}:=u^{\varepsilon}-\sum_{m=0}^{M}\varepsilon^{m}u_{m}.
\]

Relying on  the auxiliary problems \eqref{eq:prob1}-\eqref{eq:probm}, note that 
the function $\varphi^{\varepsilon}$ satisfies the following system:
\begin{equation}
\begin{cases}
\mathcal{A}^{\varepsilon}\varphi^{\varepsilon}=R\left(u^{\varepsilon}\right)-\sum_{m=0}^{M-2}\varepsilon^{m-2}R\left(u_{m}\right)\\
\quad\quad\quad\quad-\varepsilon^{M-1}\left(\mathcal{A}_{1}u_{M}+\mathcal{A}_{2}u_{M-1}\right)-\varepsilon^{M}\mathcal{A}_{2}u_{M}, & \mbox{in}\;\Omega^{\varepsilon},\\
-d^{\varepsilon}\nabla_{x}\varphi^{\varepsilon}\cdot\mbox{n}=\varepsilon^{M}d^{\varepsilon}\nabla_{x}u_{M}\cdot\mbox{n}, & \mbox{on}\;\Gamma^{\varepsilon}.
\end{cases}\label{eq:24}
\end{equation}

Now, multiplying the PDE in \eqref{eq:24} by $\varphi\in V^{\varepsilon}$
and integrating by parts, we arrive at
\begin{eqnarray}
\left\langle d^{\varepsilon}\varphi^{\varepsilon},\varphi\right\rangle _{V^{\varepsilon}} & = & \left\langle R\left(u^{\varepsilon}\right)-\sum_{m=0}^{M-2}\varepsilon^{m-2}R\left(u_{m}\right),\varphi\right\rangle _{L^{2}\left(\Omega^{\varepsilon}\right)}\nonumber \\
 &  & -\varepsilon^{M-1}\left\langle \mathcal{A}_{1}u_{M}+\mathcal{A}_{2}u_{M-1}+\varepsilon\mathcal{A}_{2}u_{M},\varphi\right\rangle _{L^{2}\left(\Omega^{\varepsilon}\right)}\nonumber \\
 &  & -\varepsilon^{M}\int_{\Omega_{int}^{\varepsilon}}d^{\varepsilon}\nabla_{x}u_{M}\cdot\mbox{n}\varphi dS_{\varepsilon}.\label{eq:4.29}
\end{eqnarray}

From here on, we have to estimate the integrals on the right-hand
side of \eqref{eq:4.29}., which is a standard procedure; see \cite{CP99,KM15} for similar calculations. Thus, we claim
that
\begin{equation}
\left|\left\langle R\left(u^{\varepsilon}\right)-\sum_{m=0}^{M-2}\varepsilon^{m-2}R\left(u_{m}\right),\varphi\right\rangle _{L^{2}\left(\Omega^{\varepsilon}\right)}\right|\le CL\left\Vert u^{\varepsilon}-\sum_{m=0}^{M}\varepsilon^{m}u_{m}+\mathcal{O}\left(\varepsilon^{M-1}\right)\right\Vert _{V^{\varepsilon}}\left\Vert \varphi\right\Vert _{L^{2}\left(\Omega^{\varepsilon}\right)},\label{eq:2.15}
\end{equation}
where we have essentially used the global Lipschitz condition on the reaction
term, the assumption \eqref{eq:newR}, and the Poincar\'e
inequality. Next, we get
\begin{equation}
\varepsilon^{M-1}\left|\left\langle \mathcal{A}_{1}u_{M}+\mathcal{A}_{2}u_{M-1}+\varepsilon\mathcal{A}_{2}u_{M},\varphi\right\rangle _{L^{2}\left(\Omega^{\varepsilon}\right)}\right|\le C\varepsilon^{M-1}\left\Vert \varphi\right\Vert _{L^{2}\left(\Omega^{\varepsilon}\right)},\label{eq:2.16}
\end{equation}
while using the trace inequality (cf. \cite[Lemma 2.31]{CP99})  to deal with the 
the last integral, it gives
\begin{equation}
\varepsilon^{M}\left|\int_{\Omega_{int}^{\varepsilon}}d^{\varepsilon}\nabla_{x}u_{M}\cdot\mbox{n}\varphi dS_{\varepsilon}\right|\le C\varepsilon^{M-1}\left\Vert \varphi\right\Vert _{L^{2}\left(\Omega^{\varepsilon}\right)}.\label{eq:2.17}
\end{equation}

Combining \eqref{eq:2.15}-\eqref{eq:2.17}, we provide that
\[
\alpha\left|\left\langle \varphi^{\varepsilon},\varphi\right\rangle _{V^{\varepsilon}}\right|\le C\varepsilon^{M-1}\left\Vert \varphi\right\Vert _{L^{2}\left(\Omega^{\varepsilon}\right)},
\]
which finally leads to $\left\Vert \varphi^{\varepsilon}\right\Vert _{V^{\varepsilon}}\le C\varepsilon^{\frac{M-1}{2}}$
by choosing $\varphi=\varphi^{\varepsilon}$, very much in the spirit of energy estimates.

Summarizing, we state our results in the frame of  the following theorems.
\begin{thm}
Suppose \eqref{eq:newR} holds for $r\in\left\{ 1,2\right\} $ and
assume $\kappa_{p}:=C_{p}L\alpha^{-1}<1$ for the given Poincar\'e
constant. Let $\left\{ u_{m}^{\left(n_{m}\right)}\right\} _{n_{m}\in\mathbb{N}}$
be the schemes that approximate the nonlinear auxiliary problems \eqref{eq:prob1-1}-\eqref{eq:probm-1}.
Then \eqref{eq:prob1-1}-\eqref{eq:probm-1} admit a unique solution
$u_{m}$ for all $m\in\left\{ 0,...,M\right\} $ with the speed of
convergence:
\[
\left\Vert u_{m}^{\left(n_{m}\right)}-u_{m}\right\Vert _{H^{1}\left(Y_{1}\right)}\le\frac{C\kappa_{p}^{n}}{1-\kappa_{p}^{n}}\;\mbox{for all}\;n_{m}\in\mathbb{N}\;\mbox{and}\;m\in\left\{ 0,...,M\right\} ,
\]
where $C>0$ is a generic $\varepsilon$-independent constant and
$n:=\max\left\{ n_{0},...,n_{M}\right\} $.
\end{thm}

\begin{thm}
Let $u^{\varepsilon}$ be the solution of the elliptic system \eqref{eq:1.1}
with the assumptions $\left(\mbox{A}_{1}\right)-\left(\mbox{A}_{2}\right)$
stated above and suppose that \eqref{eq:newR} holds for $r\in\left\{ 1,2\right\} $.
For $m\in\left\{ 0,\ldots,M\right\} $ with $M\ge2$, we consider
$u_{m}$ the solutions of the auxiliary problems \eqref{eq:prob1-1}-\eqref{eq:probm-1}.
Then we obtain the following corrector estimate:
\[
\left\Vert u^{\varepsilon}-\sum_{m=0}^{M}\varepsilon^{m}u_{m}\right\Vert _{V^{\varepsilon}}\le C\varepsilon^{\frac{M-1}{2}},
\]
where $C>0$ is a generic $\varepsilon$-independent constant.
\end{thm}

\section*{Acknowledgment}
AM thanks NWO MPE "Theoretical estimates of heat losses in geothermal wells" (grant nr. 657.014.004) for funding. VAK would like to thank Prof. Nguyen Thanh Long, his old supervisor, for the whole-hearted guidance when studying at Department of Mathematics and Computer Science, Ho Chi Minh City University of Science in Vietnam.

\bibliography{mybib}

\end{document}